\documentclass[11pt,reqno]{amsart}
\usepackage{fixltx2e}                    
\usepackage[osf]{mathpazo}  
\usepackage[utf8]{inputenc}            
\usepackage{amsmath}                     
\usepackage{amssymb, latexsym, stmaryrd, amsthm, dsfont, amsfonts, amsbsy,amsthm, amsmath, mathrsfs}            
\usepackage{mathtools}                   
\usepackage{bm}                          
\usepackage{enumerate}                   
\usepackage{verbatim}                    
\usepackage{url}   
\usepackage{lscape}                      

\usepackage{microtype}  
\usepackage[all, knot]{xy}
        \xyoption{arc} 
        \xyoption{web}                 
\makeatletter                            
\def\MT@register@subst@font{\MT@exp@one@n\MT@in@clist\font@name\MT@font@list
 \ifMT@inlist@\else\xdef\MT@font@list{\MT@font@list\font@name,}\fi}
\makeatother

\usepackage[pdftex,bookmarks,bookmarksnumbered,linktocpage,   
         colorlinks,linkcolor=blue,citecolor=blue]{hyperref}

\newcommand{\bit}{\begin{itemize}}    
\newcommand{\eit}{\end{itemize}}
\newcommand{\ben}{\begin{enumerate}}
\newcommand{\een}{\end{enumerate}}
\newcommand{\benormal}{\ben[\normalfont 1.]}   

\let\enormal\een
\newcommand{\benroman}{\ben[\normalfont (i)]}  
\let\eroman\een
\newcommand{\benbullet}{\ben[\textbullet]}     
\let\ebullet\een
\newcommand{\bde}{\begin{description}}
\newcommand{\ede}{\end{description}}

\newcommand{\?}{\ensuremath{\mkern0.4\thinmuskip}}   

\let\leq=\leqslant
\let\geq=\geqslant
\let\Box=\square                            

\let\epsilon=\varepsilon
\let\Lambda\varLambda
\let\Gamma\varGamma
\let\Delta\varDelta
\let\Lambda\varLambda
\let\Omega\varOmega
\let\Theta\varTheta
\let\Xi\varXi
\let\Pi\varPi
\let\Sigma\varSigma

\let\clsys=\mathcal                             
\let\class=\mathsf                              

\bmdefine{\A}{A}                                
\bmdefine{\2}{2}
\bmdefine{\B}{B}
\bmdefine{\D}{D}
\bmdefine{\M}{M}                                
\bmdefine{\LLL}{L}                              
\bmdefine{\Fm}{Fm}                              
\bmdefine{\zerou}{[0{,}1]}  
\bmdefine{\T}{T}                                

\newcommand{\FF}{\clsys{F}}

\newcommand{\LL}{\mathcal{L}}                   

\newcommand{\Con}{\mathrm{Con}}                            
\newcommand{\Var}{\mathnormal{V\mkern-.8\thinmuskip ar}} 
\newcommand{\FFi}{\FF\mkern-.5\thinmuskip\mathnormal{i}}  
\newcommand{\FFiLA}{\FFi_{\LL}\A}      
\newcommand{\Fi}{\text{\textsl{Fi}}}                      
\newcommand{\FiLA}{\Fi_{\LL}^{\A}}                        
\newcommand{\Mod}{\class{Mod}}

\newcommand{\Alg}{\class{Alg}}

\newcommand{\AlgL}{\Alg\?\LL}

\bmdefine{\boldstar}{\mathchoice{\textstyle*}{\textstyle*}{\textstyle*}{\scriptstyle*}}
\newcommand{\Modstar}{\Mod^{\boldstar}}
\newcommand{\ModstarL}{\Modstar\!\LL}

\newcommand{\Algstar}{\Alg^{\boldstar}\!\?}
\newcommand{\AlgstarL}{\Algstar\!\?\LL}

\bmdefine{\btau}{\tau}                                  
\bmdefine{\brho}{\rho}                                  

\newcommand{\vdashL}{\vdash_{\!\LL}}                      
\newcommand{\sineq}{\mathrel{\dashv\mkern1.5mu\vdash}}  

\bmdefine{\leibniz}{\Omega}        

\bmdefine{\frege}{\Lambda}         

\makeatletter
\newcommand{\tarskidsp}{\mathord%
   {\m@th\raisebox{0pt}[0pt][0pt]{$\stackrel%
   {\raisebox{-2.7pt}[0ex][0pt]{$\displaystyle \,\?\thicksim$}}%
   {\displaystyle\leibniz}$}}}
\newcommand{\tarskitxt}{\mathord%
   {\m@th\raisebox{0pt}[0pt][0pt]{$\stackrel%
   {\raisebox{-2.7pt}[0ex][0pt]{$\,\?\thicksim$}}{\displaystyle\leibniz}$}}}
\newcommand{\tarskiscr}{\mathord%
   {{\m@th\raisebox{0pt}[0pt][0pt]{$\stackrel%
   {\raisebox{-2.4pt}[0ex][0pt]{$\scriptstyle \,\?\thicksim$}}%
   {\scriptstyle\leibniz}$}}}}
\newcommand{\tarskiscrscr}{\mathord%
   {{\m@th\raisebox{0pt}[0pt][0pt]{$\stackrel%
   {\raisebox{-2pt}[0ex][0pt]{$\scriptscriptstyle \,\?\thicksim$}}%
   {\scriptscriptstyle\leibniz}$}}}}
\newcommand{\tarski}{\@ifnextchar ^ %
   {\mathchoice{\tarskidsp\kern-.07em}{\tarskitxt\kern-.07em}%
   {\tarskiscr\kern-.07em}{\tarskiscrscr\kern-.07em}}%
   {\mathchoice{\tarskidsp}{\tarskitxt}{\tarskiscr}{\tarskiscrscr}}}
\makeatother

\theoremstyle{theorem}
\newtheorem{Theorem}{Theorem}[section]
\newtheorem{Lemma}[Theorem]{Lemma}
\newtheorem{Corollary}[Theorem]{Corollary}

\newtheorem{claim}{\textbf{Claim}}[Theorem]

\theoremstyle{definition}
\newtheorem{Definition}[Theorem]{Definition}
\newtheorem{exa}[Theorem]{Example}

\theoremstyle{remark}

\newcommand{\lo}{\mathcal{L}} 

\newcommand{\C}{\boldsymbol{C}} 

\begin{document}
 \title{A computational glimpse at the Leibniz and Frege  hierarchies}

\author{Tommaso Moraschini}
\email{tommaso.moraschini@gmail.com}
\date{\today}

\maketitle


\begin{abstract}
In this paper we consider, from a computational point of view, the problem of classifying logics within the Leibniz and Frege hierarchies typical of abstract algebraic logic. The main result states that, for logics presented syntactically, this problem is in general undecidable. More precisely, we show that there is no algorithm that  classifies the logic of a finite consistent Hilbert calculus in the Leibniz and in the Frege hierarchies.
\end{abstract}


\begin{center}
\textit{Dedicated to Professor Josep Maria Font on the occasion of his retirement}
\end{center}

\section{Introduction}

Abstract algebraic logic (AAL for short) is a field that studies uniformly propositional logics \cite{AAL-AIT-f,Cz01,FJa09,FJaP03b}. One of its main achievements is the development of the so-called Leibniz and Frege hierarchies in which propositional logics are classified according to two different criteria. More precisely, the \textit{Leibniz hierarchy} provides a taxonomy that classifies propositional systems accordingly to the way their notions of \textit{logical equivalence} and of \textit{truth} can be defined. Roughly speaking, the location of a logic inside the Leibniz hierarchy reflects the strength of the relation that it enjoys with its algebraic counterpart. In this sense, the Leibniz hierarchy revealed to be a useful framework where to express general transfer theorems between metalogical and algebraic properties. This is the case for example for superintuitionistic logics \cite{ChZa97}, where the Beth definability property corresponds to the fact that epimorphisms are surjective  \cite{BlHoo06,BMR16,GaMa05,Ho00}, and the deductive interpolation corresponds to the amalgamation property \cite{CzP99,GaMa05,Mak77b}. On the other hand,  the \textit{Frege hierarchy} offers a classification of logics according to general \textit{replacement} principles. Remarkably, some of these replacement properties can be formulated semantically by asking that the different elements in a model of the logic are separated by a deductive filter. This is what happens for example in superintuitionistic logics, whose algebraic semantics is given by varieties of Heyting algebras where logical filters are just lattice filters.

The aim of this paper is to investigate the computational aspects of the problem of classifying syntactically presented logics in the Leibniz and Frege hierarchies. More precisely, we will consider the following problem:
\benbullet
\item Let $\class{K}$ be a level of the Leibniz (resp. Frege) hierarchy. Is it possible to decide whether the logic of a given finite consistent Hilbert calculus in a finite language belongs to $\class{K}$?
\ebullet
It turns out that in general the answer is negative both for the Leibniz and the Frege hierarchies. To show the first case, we reduce Hilbert's tenth problem on Diophantine equations to the problem of classifying the logic of a (finite consistent) Hilbert calculus in the Leibniz hierarchy, thus obtaining that also the last one is undecidable (Theorem \ref{UndecAlg}). In the process we will also describe and axiomatize a new logic, whose deductions mimic the equational theory of commutative rings with unity (Theorem \ref{Finite}). In order to prove that also the problem of classifying the logic of a (finite and consistent) Hilbert calculus in the Frege hierarchy is undecidable, we rely on the undecidability of the equational theory of relation algebras in a single variable \cite{TaGi87}. Remarkably, our proof shows that this classification problem remains undecidable even if we restrict our attention to Hilbert calculi that determine a finitary algebraizable logic (Theorem \ref{Thm:FregeUndecidable2}).

\section{Preliminaries} 

Here we present a brief survey of the main definitions and results of abstract algebraic logic we will make use of along the article; a systematic exposition can be found for example in \cite{BP86,BP89,Cz01,AAL-AIT-f,FJa09,FJaP03b}.

A \textit{closure operator} over a set $A$ is a monotone function $C \colon \mathcal{P}(A) \to \mathcal{P}(A)$ such that $X\subseteq C(X) = C\bigl(C(X)\bigr)$ for every $X\in \mathcal{P}(A)$, and a \textit{closure system} on $A$ is a family $\mathcal{C}\subseteq \mathcal{P}(A)$ 
closed under arbitrary intersections and such that $A\in \mathcal{C}$. It is well known that the closed 
sets (fixed points) of a closure operator form a closure system and that given a closure system 
$\mathcal{C}$ one can define a closure operator $C$ by letting $C(X)= \bigcap \{Y \in \mathcal{C} \?\?\vert\?\? X \subseteq Y\}$ for every $X \in \mathcal{P}(A)$. These transformations are indeed inverse to 
one another. Therefore definitions and results established for 
closure operators transfer naturally to closure systems and vice-versa. Let $\mathcal{C}\subseteq \mathcal{P}(A)$  be a closure system on $A$. We say that an element $a\in A$ is a \textit{theorem} of $\mathcal{C}$ if $a \in \bigcap \mathcal{C}$. 

Fixed an algebraic type $\mathscr{L}$ and a set $X$, we denote by $Fm(X)$ the set of formulas over $\mathscr{L}$ built up with variables in $X$ and by $\Fm(X)$ the corresponding absolutely free algebra. In particular we fix a countable set $\Var$ of variables $x, y, z$, etc., and we write $Fm$ and $\Fm$ instead of $Fm(\Var)$ and $\Fm(\Var)$. Moreover, given a formula $\varphi \in Fm$, we write $\varphi(x, \vec{z})$ if the variables of $\varphi$ are among $x$ and $\vec{z}$ and $x$ does not appear in $\vec{z}$. From now on we assume that we are working with a fixed algebraic type.

By a \textit{logic} $\LL$ we understand a closure operator $C_{\lo}\colon \mathcal{P}(Fm) \to \mathcal{P}(Fm)$ which is \textit{structural} in the sense that $\sigma C_{\lo} \subseteq C_{\lo} \sigma$ for every endomorphism (i.e. substitution) $\sigma \colon \Fm \to \Fm$. Given $\Gamma \cup \{\varphi\}\subseteq Fm$ we write $\Gamma \vdashL  \varphi$ instead of $\varphi \in C_{\lo}(\Gamma)$. Moreover, given $\Gamma \cup \{ \varphi, \psi \} \subseteq Fm$, we denote by $\Gamma, \varphi \sineq_{\LL} \psi, \Gamma$ the fact that both $\Gamma, \varphi \vdash_{\LL} \psi$ and $\Gamma, \psi \vdash_{\LL} \varphi$ are true. Given two logics $\lo$ and $\lo'$, we write $\lo \leq \lo'$ if $C_{\lo}(\Gamma) \subseteq C_{\lo'}(\Gamma)$ for every $\Gamma \subseteq Fm$. A logic $\LL$ is \textit{consistent} when there is a formula $\varphi$ such that $\emptyset \nvdash_{\LL} \varphi$. Since $\lo$ always denotes an arbitrary logic, we skip, in the formulation of our results, assumptions like ``let $\lo$ be a logic''. 

A logic $\LL$ is \textit{finitary} if for every $\Gamma \cup \{ \varphi \} \subseteq Fm$ such that $\varphi \in C_{\LL}(\Gamma)$, there is a finite subset $\Gamma' \subseteq \Gamma$ such that $\varphi \in C_{\LL}(\Gamma')$. In this paper, we will focus only on  finitary logics.  A \textit{rule} is an expression of the form $\Gamma \vdash \varphi$, where $\Gamma \cup \{ \varphi \}$ is a finite set of formulas. Given two finite sets of formulas $\Gamma$ and $\Gamma'$, the expression $\Gamma \sineq \Gamma'$ is an abbreviation for the rules $\Gamma \vdash \varphi$ and $\Gamma' \vdash \psi$ for every $\varphi \in \Gamma'$ and $\psi \in \Gamma$. A \textit{Hilbert calculus} is a set of rules. It is well known that every Hilbert calculus determines a finitary logic and, vice-versa, that every finitary logic is determined by a Hilbert calculus. An Hilbert calculus is \textit{consistent}, when so is the logic it determines. For the notion of \textit{proof} in a Hilbert calculus we refer the reader to \cite[pag. 28]{Cz01}.

We denote algebras with bold capital letters $\A$, $\B$, $\C$, etc.\ (with universes $A$, $B$, $C$, etc.\ 
respectively). We skip assumptions like ``let $\A$ be an 
algebra'' in the formulation of our results. Given an algebra $\A$, we denote its congruence lattice by $\Con\A$. The identity relation on $\A$ is denoted by $\textup{Id}_{\A}$. Given $\theta \in \Con\A$, we denote the natural surjection onto the quotient by $\pi_{\theta} \colon \A \to \A / \theta$. A congruence $\theta \in \Con\A$ is \textit{compatible} with a set $F \subseteq A$ if for every $a, b\in A$
\[
\text{if }\langle a, b \rangle \in \theta \text{ and }a \in F\text{, then }b \in F.
\]
Given $F \subseteq A$, there exists always that largest congruence on $\A$ compatible with $F$ \cite[Theorem 4.20]{AAL-AIT-f}. We 
 denote it by $\leibniz^{\A} F$ and call it the \textit{Leibniz congruence} of $F$ on 
$\A$. The Leibniz congruence can be characterized in terms of the indiscernibility of elements with respect to filters in the following way \cite[Theorem 4.23]{AAL-AIT-f}. Given an algebra $\A$, a function $p \colon A^{n} \to A$ is a \textit{polynomial function} if there are a natural number $m$, a formula $\varphi(x_{1}, \dots, x_{n+m})$, and elements $b_{1}, \dots, b_{m} \in A$ such that
\begin{align*}
p(a_{1}, \dots, a_{n}) = \varphi^{\A}(a_{1}, \dots, a_{n}, b_{1}, \dots, b_{m})
\end{align*}
for every $a_{1}, \dots, a_{n} \in A$. Observe that the notation $\varphi(x_{1},  \dots, x_{n+m})$ means just that the variables really occurring in $\varphi$ are among, but are not necessarily all, the elements of $\{x_{1}, \dots, x_{n+m}\}$. We denote by $\textup{Pol}\A$ the set of all polynomial functions of $\A$.
\begin{Lemma}\label{Polynomial} Let $\A$ be an algebra, $F \subseteq A$ and $a, b\in A$. $\langle a, b \rangle \in \leibniz^{\A}F$ if and only if for every unary $p \in \textup{Pol}\A$,
\[
p(a)\in F \Longleftrightarrow p(b)\in F.
\]
\end{Lemma}
\noindent A \textit{matrix} is a pair $\langle \A, F \rangle$ such that $F \subseteq A$. The set $F$ is called the \textit{filter} or the \textit{truth set} of the matrix. A matrix $\langle \A, F \rangle$ is \textit{reduced} if $\leibniz^{\A}F = \textup{Id}_{\A}$. 

Given a logic $\lo$ and an algebra $\A$, we 
say that a set $F \subseteq A$ is a \textit{deductive filter} of $\lo$ over $\A$ when
\begin{gather*}
\text{if }\Gamma \vdashL  \varphi\text{, then for every homomorphism }h\colon \Fm\to \A,\\ 
\text{if }h[\Gamma] \subseteq F\text{, then }h(\varphi)\in F
\end{gather*}
for every $\Gamma \cup \{\varphi\} \subseteq Fm$. We denote by $\FFiLA$ the set of deductive filters of $\lo$ over $\A$, which turns out to be a closure system \cite[Proposition 2.22]{AAL-AIT-f}. Thus, we denote by $\FiLA(\cdot)$ the closure operator of $\LL$-filter generation over the algebra $\A$. The class of  \textit{reduced models} of a given a logic $\LL$ is the following:
\[
\ModstarL \coloneqq \{ \langle \A, F \rangle : F \in \FFiLA \text{ and }\leibniz^{\A} F = \textup{Id}_{\A} \}.
\]
Given an algebra $\A$, a logic $\LL$ and a set $F \subseteq A$, the congruence
\[
\tarski^{\A}_{\LL}F \coloneqq \bigcap \{ \leibniz^{\A} G :  G  \in \FFiLA \text{ and }F \subseteq G \}
\]
is called the \textit{Suszko congruence} of $F$ over $\A$ (relative to $\LL$).  The Suszko congruence allows to associate a special class of algebras with a given logic $\LL$ as follows:
\[
\AlgL \coloneqq \{ \A : \text{ there is }F \in \FFiLA \text{ such that }\tarski_{\LL}^{\A} F = \text{Id}_{\A}\}.
\]

\section{The hierarchies of AAL}

\begin{figure}[t]
\[
\xymatrix@R=45pt @C=55pt @!0{
 &  & {\txt{finitely regurlarly \\ algebraizable }} \ar[dr]\ar[dl]  &&&\\
 & {\txt{finitely \\ algebraizable }}\ar[dr]\ar[dl] & & {\txt{regularly\\ algebraizable}}\ar[dl]\ar[dr] && \\
{\txt{finitely \\ equivalential }} \ar[dr]&& {\txt{algebraizable}}\ar[dl]\ar[dr] & & {\txt{regularly weakly\\ algebraizable}}\ar[dl]\ar[dr]& \\
&{\txt{equivalential}}\ar[dr] & & {\txt{weakly\\ algebraizable}}\ar[dl]\ar[dr] & & {\txt{assertional}}\ar[dl]\\
& & {\txt{protoalgebraic}} & & {\txt{truth-equational}}\ar[dr]&\\
&&&&&{\txt{truth is implicitly\\definable in $\ModstarL$}}
}
\]
\caption{The Leibniz hierarchy.}
\label{Hierarchy}
\end{figure}
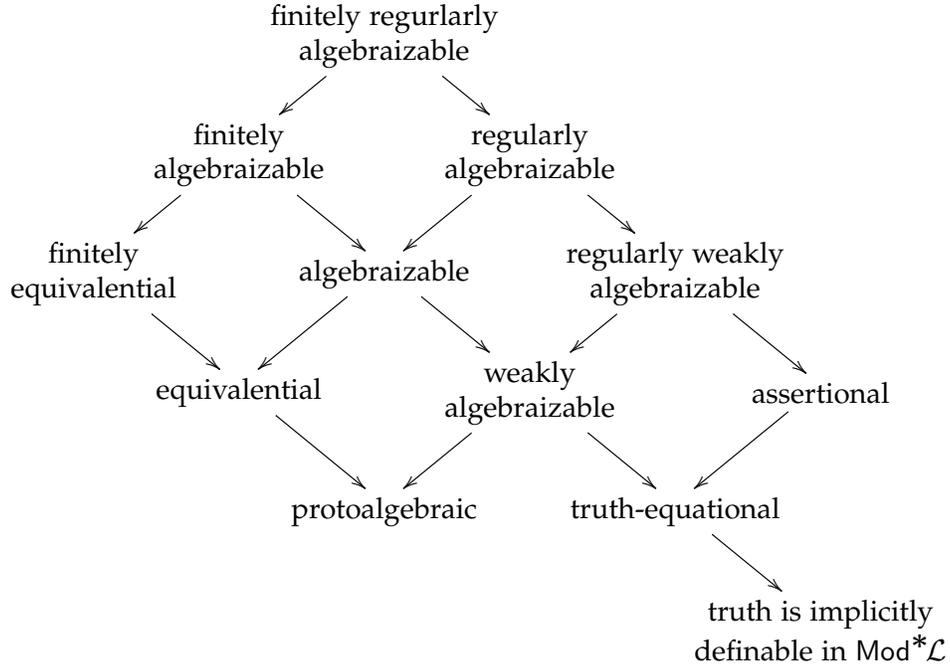

One of the main achievements of abstract algebraic logic is the development of the so-called \textit{Leibniz} 
and \textit{Frege hierarchies}, in which propositional logics $\LL$ are classified by means of properties that 
reflect the strength (or weakness) of the relation that they enjoy with $\AlgL$. We 
 begin by describing the main classes of logics that constitute the Leibniz hierarchy (see Figure \ref{Hierarchy}, where arrows denote the inclusion relation).   
In order to 
 follow a uniform approach we chose to do this from a semantic point of view, even if some of these 
 concepts were originally defined in some different ways. The main idea, which lies behind the 
 development of the Leibniz hierarchy, is that propositional logics can be characterized by means of two 
 basic features: the fact that logical equivalence (which is identified with the Leibniz congruence) 
 is definable by means of formulas, and the fact that the truth sets in $\ModstarL$ are definable either implicitly or by 
 means of equations. 

A logic $\LL$ is \textit{protoalgebraic} if there is a set of formulas $\brho(x, y, \vec{z})$ in two variables $x$ and $y$ (possibly) with parameters $\vec{z}$ such that for every algebra $\A$, every $F \in \FFiLA$ and every $a, b \in A$
\begin{equation}\label{Defines}
\langle a, b  \rangle \in \leibniz^{\A}F\Longleftrightarrow\brho^{\A}(a, b, \vec{c}) \subseteq F\text{ for every }\vec{c} \in A.
\end{equation}
In this case $\brho(x, y, \vec{z})$ is a set of \textit{congruence formulas with parameters} for $\LL$. Analogously, a logic is called \textit{equivalential} if there is set of  formulas $\brho(x, y)$ in two variables and without parameters which satisfies (\ref{Defines}). In this case $\brho(x, y, \vec{z})$ is called a set of \textit{congruence formulas} for $\LL$. Protoalgebraic logics can be characterized syntactically as follows \cite[Proposition 6.7 and Theorem 6.57]{AAL-AIT-f}. A set of formulas in two variables $\brho(x, y)$ is a set of \textit{protoimplication formulas} for a logic $\LL$ if the following conditions hold:
\begin{align*}
 \emptyset &\vdash_{\LL} \brho(x, x)\tag{\textup{R}} \\
 x, \brho(x, y) &\vdash_{\LL} y \tag{\textup{MP}}
\end{align*} 

\begin{Theorem}\label{SynProt}
A logic $\LL$ is protoalgebraic if and only if it has a set of protoimplication formulas.
\end{Theorem}

\textit{Truth is equationally definable} in a class of matrices $\class{M}$ if there is a set $\btau(x)$ of equations in variable $x$ such that for every $\langle \A, F \rangle \in \class{M}$,
\[
F = \{ a \in A : \A \vDash \btau(a) \}.
\]
\textit{Truth is implicitly definable} in a class of matrices $\class{M}$ if the matrices in $\class{M}$ are determined by their algebraic reducts, in the sense that if $\langle \A, F\rangle, \langle \A, G \rangle \in \class{M}$, then $F = G$. Clearly equational definability implies implicit definability. A logic $\LL$ is 
\textit{truth-equational} if truth is equationally definable in $\ModstarL$.

A logic $\LL$ is \textit{weakly algebraizable} if it is both protoalgebraic and truth-equational. It is \textit{algebraizable} if it is both equivalential and truth-equational. Algebraizable logics enjoy the following syntactic characterization \cite[Theorem 3.21]{AAL-AIT-f}:

\begin{Theorem}\label{SynAlg} A logic $\LL$ is algebraizable if and only if there are a set of equations $\btau(x)$ and a set of formulas $\brho(x, y)$ such that:
\begin{align*}
\emptyset \vdash_{\LL}& \brho(x, x) \tag{\textup{R}}\\
x, \brho(x, y)\vdash_{\LL}& y \tag{\textup{MP}}\\
\bigcup_{i \leq n}\brho(x_{i}, y_{i})\vdash_{\LL}& \brho(\lambda x_{1}\dots x_{n}, \lambda y_{1} \dots y_{n}) \tag{\textup{Rep}}\\
x \sineq_{\LL}& \{ \varphi(\epsilon, \delta) : \varphi(x, y)\in \brho\text{ and }\epsilon \thickapprox \delta \in \btau \} \tag{\textup{A3}}
\end{align*}
for every $n$-ary function symbol $\lambda$. In this case $\brho(x, y)$ a set of congruence formulas for $\LL$.
\end{Theorem}
\begin{Lemma}\label{Lem:MovingEqIntoFm}
Let $\LL$ be an algebraizable logic with set of congruence formulas $\brho(x, y)$. For every $\varphi, \psi \in Fm$, we have
\[
\emptyset \vdash_{\LL}\brho(\varphi, \psi) \Longleftrightarrow \AlgL \vDash \varphi \thickapprox \psi.
\]
\end{Lemma}

A logic $\LL$ is \textit{assertional} if there is a class $\class{K}$ of algebras with a constant term $1$ such that for every $\Gamma \cup \{ \varphi \} \subseteq Fm$,
\begin{align*}
\Gamma \vdash_{\LL} \varphi \Longleftrightarrow& \text{ for every }\A \in \class{K}\text{ and homomorphism }h \colon \Fm \to \A,\\
&\text{ if }h[\Gamma] \subseteq \{ 1 \} \text{, then }h(\varphi) = 1.
\end{align*}
\noindent A logic $\LL$ is \textit{regularly }(\textit{weakly}) \textit{algebraizable} if it is assertional and (weakly) algebraizable. Regularly algebraizable logics enjoy a syntactic characterization \cite[Theorem 3.52]{AAL-AIT-f}.

\begin{Lemma}\label{AssertionalDelta}
Let $\LL$ be an algebraizable logic. $\LL$ is regularly algebraizable if and only if
\begin{align*}
x, y &\vdash_{\LL}\brho(x, y) \tag{\textup{G}}
\end{align*}
for some (or, equivalently, every) of its sets of congruence formulas $\brho(x, y)$.
\end{Lemma}

\noindent Finally a logic $\LL$ is \textit{finitely equivalential }(resp. \textit{finitely algebraizable, finitely regularly algebraizable})  if it is equivalential (resp. algebraizable, regularly algebraizable) and it has a \textit{finite} set of congruence formulas.

\begin{figure}[t]
\[
\xymatrix@R=47pt @C=57pt @!0{
 & {\txt{fully Fregean }} \ar[dr]\ar[dl]  & \\
{\txt{fully \\ selfextensional }} \ar[dr] & & {\txt{Fregean }} \ar[dl] \\
&{\txt{selfextensional }}&
}
\]
\caption{The Frege hierarchy.}
\label{Frege hierarchy fig}
\end{figure}
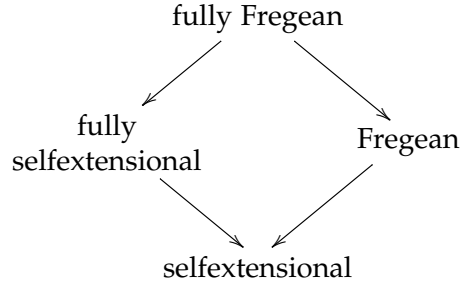

Now we turn to describe briefly the structure of the other hierarchy typical of abstract algebraic logic, namely the Frege hierarchy (see Figure \ref{Frege hierarchy fig}, where arrows denote the inclusion relation).
In this hierarchy logics are classified by means of general replacement properties, but we choose here to define them by the most convenient equivalent characterization. A logic $\LL$ is \textit{selfextensional} if the relation $\sineq_{\LL}$ is a congruence on $\Fm$. In particular, we will make use of the following easy observation \cite[Proposition 2.43]{FJa09}. 
\begin{Lemma}\label{EquationalFrege}
Let $\LL$ be selfextensional and $\alpha, \beta \in Fm$. If $\alpha \sineq_{\LL}\beta$, then $\AlgL \vDash \alpha \thickapprox \beta$. 
\end{Lemma}
\noindent A logic $\LL$ is \textit{fully selfextensional} if for every $\A \in \AlgL$ and $a, b \in A$:
\[
\text{if }\FiLA\{a \} =\FiLA\{ b \}\text{, then }a=b.
\]
Fully selfextensional logics are also selfextensional, while the converse is not true in general \cite{Ba03}. A logic $\LL$ is \textit{Fregean} if for every $n$-ary function symbol $\lambda$ and formulas $\Gamma \cup \{ \alpha_{1}, \dots, \alpha_{n}, \beta_{1}, \dots, \beta_{n} \} \subseteq Fm$, if $\Gamma, \alpha_{i} \sineq_{\LL} \beta_{i}, \Gamma$ for every $i \leq n$, then $\Gamma, \lambda(\alpha_{1}, \dots, \alpha_{n}) \sineq_{\LL} \lambda( \beta_{1}, \dots, \beta_{n}), \Gamma$. Finally $\LL$ is \textit{fully Fregean} if for every algebra $\A \in \AlgL$ and every $F \in \FFiLA$ we have that
\[
\text{if }\FiLA(F \cup \{ a \}) = \FiLA(F \cup \{ b \})\text{, then }\langle a, b \rangle \in \tarski_{\LL}^{\A}F.
\]

\section{The classification problem in the Leibniz hierarchy}

We will begin our tour around some computational aspects of abstract algebraic logic by considering the problem of classifying logics in the Leibniz hierarchy. The goal of this part is to prove that the problem of classifying logics determined by a finite consistent Hilbert calculus in a finite language in the Leibniz hierarchy is in general undecidable (Theorem \ref{UndecAlg}). Our strategy is the following: for every Diophantine equation $p\thickapprox 0$ we define a finite consistent Hilbert calculus $\LL(p)$ such that $\LL(p)$ belongs to a given level of the Leibniz hierarchy if and only if $p\thickapprox 0$ has an integer solution (Definition \ref{Def:Logic} and Theorem \ref{Thm:LeibnizUndecidable}). For the sake of completeness we recall some concepts. An algebra $\A = \langle A, +, \cdot, -, 0, 1 \rangle$ of type $\langle 2, 2, 1, 0, 0 \rangle$ is a \textit{commutative ring} if $\langle A, +, -, 0\rangle$ is an Abelian group, $\langle A, \cdot, 1 \rangle$ is a commutative monoid, and
\[
a \cdot ( b + c ) = ( a \cdot b ) + ( a \cdot c )  
\]
for every $a, b, c \in A$. We denote by $\class{CR}$ the variety of commutative rings. The algebra $\mathds{Z} \coloneqq \langle Z, +, \cdot, -, 0, 1 \rangle$, where $Z$ is the collection of integer numbers, is a commutative ring. A \textit{Diophantine equation} is an equation of the form $p(z_{1}, \dots, z_{1}) \thickapprox 0$, where $p(z_{1}, \dots, z_{1})$ is a term in the language of commutative rings. Hilbert's tenth problem asked for an algorithm that, given a Diophantine equation, tells us whether it has a solution in $\mathds{Z}$ or not. Matiyasevich showed that such an algorithm does not exist \cite{Mati93}, see also  \cite{BeMa77}. In other words, the problem of determining whether a given Diophantine equation has an integer solution is undecidable.

In order to relate this problem to the one of the classification of logics into the Leibniz hierarchy, we will construct a logic that mimics the behaviour of (Diophantine) equations  in commutative rings. In \cite{FMo14c} (see also \cite{Mo13a}) a way of doing this for arbitrary varieties is described. More precisely, given a non-trivial variety $\class{V}$, we let $\LL_{\class{V}}$ be  the logic determined by the following class of matrices:
\[
\{ \langle \A, F \rangle : \A \in \class{V}\text{ and }F \subseteq A \}.
\]
Given $\Gamma \cup \{ \varphi \} \subseteq Fm$, we will write $\Gamma \vdash_{\class{V}} \varphi$ as a shortening for $\Gamma \vdash_{\LL_{\class{V}}}\varphi$. The following result will be used later on:
\begin{Lemma}\label{DeductionLemma}
Let $\class{V}$ be a non-trivial variety and $\Gamma \cup \{ \varphi \} \subseteq Fm$.
\benormal
\item $\AlgstarL_{\class{V}} \subseteq \AlgL_{\class{V}} = \class{V}$.
\item $\LL_{\class{V}}$ is fully selfextensional.
\item $\Gamma \vdash_{\class{V}}\varphi$ if and only if there is $\gamma\in \Gamma$ such that $\class{V} \vDash \gamma \thickapprox \varphi$.
\enormal
\end{Lemma}
\begin{proof}
Point 1 and 2 are respectively Theorem 3.5 and Corollary 3.6 of \cite{FMo14c}. Point 3, even if not explicit stated there, is a direct consequence of Lemmas 3.3 and 3.9 of \cite{FMo14c}.
\end{proof}

\noindent \textbf{Digression.} Next we are going to axiomatize the logic $\LL_{\class{CR}}$, but first let us consider the general problem of axiomatizing the logic $\LL_{\class{V}}$ for an arbitrary non-trivial variety $\class{V}$.  Suppose that we are given an equational basis $\Sigma$ for $\class{V}$. It would be nice to have a natural way of transforming $\Sigma$ into an axiomatization of the logic $\LL_{\class{V}}$. This  can be done easily, adapting Birkhoff's Completeness Theorem of equational logic \cite[Theorem 14.19, Section II]{BuSa00}, if what we are looking for is a Gentzen system adequate to $\LL_{\class{V}}$ \cite[Theorem~1.2]{FMo14d}.\footnote{Even if we won't pursue this here, it is possible to show that there is no Gentzen system \textit{fully adequate} to $\LL_{\class{V}}$ in the sense of \cite{FJa09}.} On the contrary, it is not straightforward to build a nice Hilbert calculus for $\LL_{\class{V}}$ out of $\Sigma$. In particular, the obvious idea of considering the logic axiomatized by the rules $\alpha \sineq \beta$ for all $\alpha \thickapprox \beta \in \Sigma$ does not work in general. For example, it is the case that the logic $\mathcal{S}$ determined by the rules
\begin{equation}\label{SemiLat}
x \sineq x \land x \qquad x\land y \sineq y \land x \qquad x \land ( y \land z ) \sineq (x \land y) \land z
\end{equation}
is not the logic $\LL_{\class{SL}}$ associated with the variety of semilattices $\class{SL}$ (cfr. \cite[Example~3.8]{FMo14c}). This is because the matrix $\langle \mathds{Z}_{3}, \{1, 2 \} \rangle$, where $\mathds{Z}_{3}$ is the additive semigroup of integers modulo 3, is easily proved to be a reduced model of $\mathcal{S}$. In particular this implies that $\mathds{Z}_{3} \in \Alg\mathcal{S}$ and, therefore, that $\Alg\mathcal{S} \ne \class{SL}$. Hence $\mathcal{S} \ne \LL_{\class{SL}}$, by point 1 of Lemma \ref{DeductionLemma}. In order to obtain a complete Hilbert-style axiomatization of $\LL_{\class{SL}}$ one has to add to the rules in (\ref{SemiLat}) the following ones:
\begin{align*}
u \land x &\sineq u \land (x \land x) \\
 u \land (x\land y) &\sineq u \land (y \land x) \\
  u \land (x \land ( y \land z )) &\sineq u \land ((x \land y) \land z)
\end{align*}
This is mainly due to the fact that selfextensionality, which is not easily expressible by means of Hilbert-style rules, fails for $\mathcal{S}$, while $\LL_{\class{SL}}$ is selfextensional by point 2 of Lemma \ref{DeductionLemma}.

Even if it is not straightforward to present an explicit Hilbert-style axiomatization of $\LL_{\class{V}}$, given a base $\Sigma$ for $\class{V}$, one may wonder whether there exists (no matter which one) a finite Hilbert-style axiomatization of $\LL_{\class{V}}$ when $\Sigma$ is finite. In the next example we show that in general this is not the case.

\begin{exa}[\textsf{Finite Axiomatizability}]\label{Magma}
An algebra $\A = \langle A, \cdot \rangle$ is a \textit{commutative magma} if $\cdot$ is a binary commutative operation. Clearly the class of commutative magmas forms a finitely based variety, which we denote by $\class{CM}$. We will prove that the logic $\LL_{\class{CM}}$ is not axiomatizable by means of a finite set of Hilbert-style rules. In order to do this, let $\mathcal{CM}$ be a finite set of rules holding in $\LL_{\class{CM}}$. We will show that there is a model of $\mathcal{CM}$ that is not a model of $\LL_{\class{CM}}$. First observe that there is a natural number $n$ that bounds the number of occurrences of (possibly equal) variables in terms appearing in the rules of $\mathcal{CM}$. We can assume, without loss of generality, that $n \geq 2$. Then we consider the algebra $\A = \langle \{0, 1, 2, \dots, n \}, \cdot \rangle$ equipped with a binary operation such that $1 \cdot 2 \coloneqq 2$ and $2 \cdot 1 \coloneqq 1$ and
\begin{displaymath}
a \cdot b = b \cdot a \coloneqq \left\{\begin{array}{@{\,}ll}
a & \text{if $a \ne n$ and $b = 0$}\\
0 & \text{if $a = n$ and $b = 0$}\\
a & \text{if $b = a - 1$ and $a \geq 3$}\\
a - 1 & \text{if $b = a - 2$ and $a \geq 3$}\\
1 & \text{otherwise}\\
\end{array} \right.
\end{displaymath}
for every $a, b \in A$ such that $\{a, b \} \ne \{1, 2 \}$.

We first show that $\langle \A, \{ 0 \} \rangle$ is not a model of $\LL_{\class{CM}}$. Observe that $\A \notin \class{CM}$, since $1\cdot 2 \ne 2 \cdot 1$. By point 1 of Lemma \ref{DeductionLemma} we know that it will be enough to prove that $\langle \A, \{0 \}\rangle$ is a reduced matrix. By point 1 of Lemma \ref{Polynomial} this amounts to checking whether for every different $a, b\in A \smallsetminus \{ 0 \}$ there is a polynomial function $p\colon \A \to \A$ such that $p(a)= 0$ if and only if $p(b) \ne 0$. This is what we do now: consider a pair of different $a, b\in A \smallsetminus \{ 0 \}$. Assume, without loss of generality, that $a < b$. Then we consider the polynomial function
\begin{align*}
p(x) \coloneqq& (\dots(((\dots((\dots((1 \cdot 2) \cdot 3) \cdot \dots ) \cdot a) \cdot \ldots \\
&\ldots\cdot (b-1)) \cdot x)\cdot (b+1)) \cdot \ldots \cdot n ) \cdot 0.
\end{align*}
It is easy to see that $p(b) = 0$. Then we turn to show that $p(a) \ne 0$. We consider two cases, whether $b-1 < 3$ or not. First consider the case in which $b-1 < 3$. We have that either ($a= 1$ and $b=2$) or ($a= 1$ and $b=3$) or ($a= 2$ and $b=3$). It is easy to prove that
\begin{align*}
\text{if }a=1\text{ and }b=2&\text{, then }p(a)=n-1\\
\text{if }a=1\text{ and }b=3&\text{, then }p(a)=1\\
\text{if }a=2\text{ and }b=3&\text{, then }p(a)=1.
\end{align*}
Then we turn to the case in which $3 \leq b-1$. We have that:
\begin{displaymath}
p(a) = (\dots((b-1 \cdot a) \cdot b+1)\ldots \cdot n) \cdot 0 =  \left\{\begin{array}{@{\,}ll}
n-1 & \text{if $a = b-2 $}\\
1 & \text{otherwise.}\\
\end{array} \right.
\end{displaymath}
Therefore we obtain that $p(a)\ne 0$. This concludes the proof that $\langle \A, \{ 0 \} \rangle$ is not a model of $\LL_{\class{CM}}$.

Now we turn to prove that $\langle \A, \{ 0 \} \rangle$ is a model of $\mathcal{CM}$. Consider a rule $\Gamma \vdash \varphi$ in $\mathcal{CM}$. Suppose that there is a homomorphism $h \colon \Fm \to \A$ such that $h[\Gamma] \subseteq \{ 0 \}$. From point 3 of Lemma \ref{DeductionLemma} we know that there is $\gamma \in \Gamma$ such that $\class{CM} \vDash \gamma \thickapprox \varphi$. In particular, we have that $h(\gamma) = 0$. We claim that $\{ h(\epsilon), h(\delta)\} \ne \{1, 2 \}$ for every subformula $\epsilon \cdot \delta$ of $\gamma$. Suppose the contrary towards a contradiction. Then there is a subformula $\epsilon \cdot \delta$ of $\gamma$ such that $\{ h(\epsilon), h(\delta)\} = \{1, 2 \}$. Since at most $n$ (possibly equal) variables occur in $\gamma$, and $\epsilon \cdot \delta$ contains at least two of them, we know that if we draw the subformulas tree of $\gamma$ there are at most $n-2$ nodes $\gamma_{1} < \gamma_{2} < \dots < \gamma_{n-2}$ (with $\gamma_{n-2} = \gamma$) strictly above $\epsilon \cdot \delta$, where $<$ is the strict order of the subformulas tree. From the definition of $\cdot$ it follows that $1 \leq h(\gamma_{m}) \leq m+2$ for every $m \leq n-2$. Therefore we obtain that $1 \leq h(\gamma) \leq n$, against the assumption that $h(\gamma) = 0$. This concludes the proof of our claim.

Now, recall that $\class{CM} \vDash \gamma \thickapprox \varphi$. It is possible to prove by induction on the length of the proofs of Birkhoff's equational logic that $\varphi$ is obtained from $\gamma$ in the following way. We replace a subformula $\gamma_{1}\cdot \gamma_{2}$ of $\gamma$ by $\gamma_{2}\cdot \gamma_{1}$ and denote by $\gamma'$ the formula obtained in this way. Then we repeat this process on $\gamma'$. Iterating this construction a finite number of times we reach $\varphi$. Now, observe that the operation $\cdot$ in $\A$ commutes always except for the case in which its arguments exhaust the set $\{1, 2 \}$. This fact, together with our claim and the observation on equational logic, implies that $h(\gamma) = h(\varphi)$. In particular, this means that $h(\varphi) = 0$. Hence $\langle \A, \{ 0 \} \rangle$ is a model of the rules in $\mathcal{CM}$.
\qed
\end{exa}
\vspace{-0.4cm}
\begin{flushright}
\textbf{End of digression}
\end{flushright}

The digression shows that the quest for a finite Hilbert-style axiomatization of the logic $\LL_{\class{CR}}$ is not in principle a trivial one: such an axiomatization may even fail to exist, as in the case of commutative magmas of Example \ref{Magma}. Nevertheless we will provide an explicit and finite Hilbert calculus for the logic $\LL_{\class{CR}}$ (Theorem \ref{Finite}).
\begin{Definition}
Let $\mathcal{CR}$ be the following Hilbert calculus (and the logic it determines) in the language of commutative rings:

\begin{align*}
w + (u \cdot ( (x \cdot y) \cdot z) ) &\sineq  w + (u \cdot ( x \cdot (y \cdot z) ) \tag{\textup{A}}\\
w + (u \cdot ( x \cdot y)) &\sineq  w + (u \cdot ( y \cdot x))  \tag{\textup{B}}\\
w + (u \cdot ( x \cdot 1)) &\sineq  w + (u \cdot x) \tag{\textup{C}}\\
w + (u \cdot ( (x + y) + z)) &\sineq w + (u \cdot ( x + ( y + z ))) \tag{\textup{D}}\\
w + (u \cdot ( x + y )) &\sineq  w + (u \cdot ( y + x )) \tag{\textup{E}}\\
w + (u \cdot ( x + 0 )) &\sineq w + (u \cdot x ) \tag{\textup{F}}\\
w + (u \cdot ( x + -  x)) &\sineq  w + (u \cdot 0 ) \tag{\textup{G}}\\
w + (u \cdot ( x \cdot ( y + z ))) &\sineq w + (u \cdot ( ( x \cdot y ) + (x \cdot z ) )) \tag{\textup{H}}\\
w + (u \cdot - ( x + y)) &\sineq w + (u \cdot ( - x + - y )) \tag{\textup{I}}\\
w + (u \cdot - ( x \cdot y)) &\sineq w + (u \cdot ( - x \cdot y )) \tag{\textup{L}}\\
w + (u \cdot - ( x \cdot y)) &\sineq w + (u \cdot (  x \cdot - y )) \tag{\textup{M}}\\
0 + x  &\sineq x \tag{\textup{N}}\\
 x + ( 1 \cdot y)   &\sineq x + y  \tag{\textup{O}}\\
\end{align*}
\end{Definition}

 The following technical result about the special nature of deductions of $\mathcal{CR}$, will be useful in the subsequent proofs.
\begin{Lemma}\label{FiniteProofs} Let $\Gamma \cup \{ \varphi \} \subseteq Fm$. If $\Gamma \vdash_{\mathcal{CR}}\varphi$, then there is $\gamma \in \Gamma$ and a finite sequence of formulas $\langle \alpha_{1}, \alpha_{2}, \dots, \alpha_{n}\rangle$ such that $\alpha_{1} = \gamma$, $\alpha_{n} = \varphi$ and for every $m< n$ there is a rule $\epsilon \sineq \delta$ of the calculus $\mathcal{CR}$ and a substitution $\sigma$ such that $\{\sigma\epsilon, \sigma \delta\} = \{\alpha_{m}, \alpha_{m+1}\}$.
\end{Lemma}
\begin{proof}
The result follows from two observations. First, that the calculus $\mathcal{CR}$ consists of rules with only one premise. Second, that if $\alpha \vdash \beta \in \mathcal{CR}$, then $\beta \vdash \alpha \in \mathcal{CR}$ too.
\end{proof}

\noindent Our main goal will be to prove that $\mathcal{CR}$ is in fact a finite axiomatization of $\LL_{\class{CR}}$. To do this we will make use of the following lemma, whose easy proof (contained in the Appendix) involves some tedious calculations. The reader may safely choose to skip it in order keep track of the proof of the main result of the section.

\begin{Lemma}\label{CR is selfextensional}
The logic $\mathcal{CR}$ is selfextensional.
\end{Lemma}

Drawing consequences from the fact that $\mathcal{CR}$ is selfextensional, the following result is easy to obtain.

\begin{Theorem}\label{Finite} The rules $\mathcal{CR}$ provide a finite axiomatization of $\LL_{\class{CR}}$.
\end{Theorem}
\begin{proof}
First observe that each of the rules in $\mathcal{CR}$ corresponds to an equation, which holds in $\class{CR}$. Together with point 3 of Lemma \ref{DeductionLemma}, this implies that $\mathcal{CR} \leq \LL_{\class{CR}}$. Now, applying (N) and (O), it is easy to prove that
\begin{align*}
(x \cdot y) \cdot z  &\sineq_{\mathcal{CR}}  x \cdot (y \cdot z) \\
x \cdot y&\sineq_{\mathcal{CR}} y \cdot x \\
x \cdot 1 &\sineq_{\mathcal{CR}}  x\\
(x + y) + z &\sineq_{\mathcal{CR}} x + ( y + z )\\
x + y  &\sineq_{\mathcal{CR}}  y + x \\
x + 0 &\sineq_{\mathcal{CR}} x \\
x + -  x &\sineq_{\mathcal{CR}} 0 \\
 x \cdot ( y + z ) &\sineq_{\mathcal{CR}}  ( x \cdot y ) + (x \cdot z ) 
\end{align*}
Since $\mathcal{CR}$ is selfextensional by Lemma \ref{CR is selfextensional}, we can apply Lemma \ref{EquationalFrege} obtaining $\Alg\mathcal{CR} \subseteq \class{CR}$. On the other hand, from point 1 of Lemma \ref{DeductionLemma} and the fact that $\mathcal{CR} \leq \LL_{\class{CR}}$, it follows that $\class{CR} = \AlgL_{\class{CR}} \subseteq \Alg\mathcal{CR}$. Therefore we conclude that $\Alg\mathcal{CR} = \class{CR}$. Since for every $\A \in \class{CR}$ and $F \subseteq A$ the matrix $\langle \A, F \rangle$ is a model of $\LL_{\class{CR}}$, this implies that $\LL_{\class{CR}} \leq \mathcal{CR}$. Thus we conclude that $\mathcal{CR} = \LL_{\class{CR}}$ as desired.
\end{proof}

Now that we have built the machinery necessary to speak about commutative rings by means of the propositional logic $\LL_{\class{CR}}$, we turn back to the classification of logics determined by finite  Hilbert calculi in the Leibniz hierarchy. We begin by describing a general way of associating a logic with each Diophantine equation.

\begin{Definition}\label{Def:Logic}
Let $p(z_{1}, \dots, z_{n}) \thickapprox 0$ be a Diophantine equation and $x$ and $y$ two new different variables. Let $\LL(\?p)$ be the logic, in the language of commutative rings expanded with a new binary operation symbol $\leftrightarrow$, determined by the following Hilbert calculus:
\begin{align*}
\emptyset \vdash & x \leftrightarrow x \tag{\textup{R}}\\
x \leftrightarrow y \vdash &y \leftrightarrow x \tag{\textup{S}}\\
x \leftrightarrow y, y \leftrightarrow z \vdash & x  \leftrightarrow z \tag{\textup{T}}\\
x \leftrightarrow y\vdash & -x \leftrightarrow -y \tag{\textup{Rep1}}\\
x \leftrightarrow y, z \leftrightarrow u \vdash &(x+z)\leftrightarrow (y + u) \tag{\textup{Rep2}}\\
x \leftrightarrow y, z \leftrightarrow u \vdash &(x \cdot z)\leftrightarrow (y \cdot u) \tag{\textup{Rep3}}\\
x \leftrightarrow y, z \leftrightarrow u \vdash &(x \leftrightarrow z)\leftrightarrow (y \leftrightarrow u) \tag{\textup{Rep4}}\\
p(z_{1}, \dots, z_{n}) \leftrightarrow 0, x, x \leftrightarrow y  \vdash & y \tag{\textup{MP'}}\\
p(z_{1}, \dots, z_{n})\leftrightarrow 0, x \sineq & x \leftrightarrow ( x \leftrightarrow x), p(z_{1}, \dots, z_{n})\leftrightarrow 0 \tag{\textup{A3'}}\\
p(z_{1}, \dots, z_{n}) \leftrightarrow 0, x, y \vdash& x \leftrightarrow y \tag{\textup{G'}}\\
\emptyset \vdash& \alpha \leftrightarrow \beta \tag{\textup{CR}}
\end{align*}
for every $\alpha \sineq \beta 	\in \mathcal{CR}$.
\end{Definition}
Observe that $\LL(\?p)$ is determined by an explicit finite Hilbert calculus in a finite language. It turns out that there is a strong relation between the existence of an integer solution to $p(z_{1}, \dots, z_{n}) \thickapprox 0$ and the location of $\LL(\? p)$ in the Leibniz hierarchy.

\begin{Theorem}\label{Thm:LeibnizUndecidable}
Let $p(z_{1}, \dots, z_{n}) \thickapprox 0$ be a Diophantine equation. The following conditions are equivalent:

\benroman
\item $\LL(\? p )$ is finitely regularly algebraizable.
\item Truth is implicitly definable in $\Modstar\LL(\? p )$.
\item $\LL(\? p )$ is protoalgebraic.
\item The equation $p(z_{1}, \dots, z_{n}) \thickapprox 0$ has an integer solution.
\eroman
\end{Theorem}
\begin{proof}
Clearly (i) implies (ii) and (iii). Now we turn to prove (ii)$\Rightarrow$(iv). We reason by contraposition. Suppose that $p(z_{1}, \dots, z_{n}) \thickapprox 0$ has no integer solution. Now choose two different integers $s$ and $m$. Then let $\boldsymbol{Z}$ be the expansion of $\mathds{Z}$ with a new binary function $\leftrightarrow$ defined as follows:
\begin{displaymath}
a \leftrightarrow b \coloneqq \left\{\begin{array}{@{\,}ll}
s & \text{if $a = b$}\\
m & \text{otherwise}\\
\end{array} \right.
\end{displaymath}
for every $a, b \in Z$. Pick $k \notin \{ s, m\}$. It is easy to check that $\langle \boldsymbol{Z}, \{ s \} \rangle$ and $\langle \boldsymbol{Z}, \{ s, k \} \rangle$ are models of $\LL(\? p )$, since $p(z_{1}, \dots, z_{n}) \thickapprox 0$ has no integer solution. Moreover they are reduced. In order to see this, pick two different $a, b \in Z$ and consider the polynomial function $q(z)\coloneqq a \leftrightarrow z$. We have that $q(a) = s$ and that $q(b) = m \notin \{ s, k \}$. By point 1 of Lemma \ref{Polynomial} we conclude both that $\langle a, b \rangle \notin \leibniz^{\boldsymbol{Z}}\{s \}$ and that $\langle a, b \rangle \notin \leibniz^{\boldsymbol{Z}}\{s, k \}$. But this means that there are two different reduced models of $\LL(\? p )$ with the same algebraic reduct and, therefore, that truth is not implicitly definable in $\ModstarL(\? p )$.

Now we prove (iii)$\Rightarrow$(iv). Suppose that $\LL(\? p )$ is protoalgebraic. Then there is a set of protoimplication formulas $\brho(x, y)$. In particular, we have that $x, \brho(x, y) \vdash_{\LL(\? p )} y$. Then there is a finite proof $\pi$ of $y$ from the premises in $\{x \} \cup \brho(x, y)$. Taking a closer look at the axiomatization of $\LL(\? p )$, it is easy to see that either an application of (MP') or of (A3') must occur in $\pi$. This is because the other rules yield complex conclusions. In particular, this implies that there is a substitution $\sigma$ such that $x, \brho(x, y) \vdash_{\LL(\? p )} \sigma p(z_{1}, \dots, z_{n}) \leftrightarrow \sigma 0$. From Theorem \ref{SynProt} we know that $\emptyset \vdash_{\LL(\? p )} \brho(x, x)$. In particular, this means that $x \vdash_{\LL(\? p )} \sigma_{x}\sigma p(z_{1}, \dots, z_{n}) \leftrightarrow \sigma_{x}\sigma 0$, where $\sigma_{x}$ is the substitution which sends all variables to $x$. But observe that every substitution leaves $0$ fixed. Therefore we conclude that 
\begin{equation}\label{MainObservation}
x \vdash_{\LL(\? p )} \sigma_{x}\sigma p(z_{1}, \dots, z_{n}) \leftrightarrow  0.
\end{equation}
Now we consider the algebra $\boldsymbol{Z}$ built in the proof of part (ii)$\Rightarrow$(iv). It is easy to check that $\langle \boldsymbol{Z}, \{ s \} \rangle$ is a model of $\LL(\? p )$. Therefore, pick an homomorphism $h\colon \Fm \to \boldsymbol{Z}$ that sends $x$ to $s$. From (\ref{MainObservation}) it follows that $h(\sigma_{x}\sigma p(z_{1}, \dots, z_{n}) \leftrightarrow  0) = s$. But observe that
\begin{align*}
h(\sigma_{x}\sigma p(z_{1}, \dots, z_{n}) \leftrightarrow  0) = s &\Longleftrightarrow h(p(\sigma_{x}\sigma(z_{1}), \dots, \sigma_{x}\sigma(z_{n})) \leftrightarrow  0) = s\\
&\Longleftrightarrow p^{\boldsymbol{Z}}(h\sigma_{x}\sigma(z_{1}), \dots, h\sigma_{x}\sigma(z_{n})) \leftrightarrow  h0 = s\\
&\Longleftrightarrow p^{\boldsymbol{Z}}(h\sigma_{x}\sigma(z_{1}), \dots, h\sigma_{x}\sigma(z_{n})) =  h0\\
&\Longleftrightarrow p^{\boldsymbol{Z}}(h\sigma_{x}\sigma(z_{1}), \dots, h\sigma_{x}\sigma(z_{n})) =  0.
\end{align*}
Therefore we conclude that $\langle h\sigma_{x}\sigma(z_{1}), \dots, h\sigma_{x}\sigma(z_{n})\rangle$ is an integer solution to the equation $p(z_{1}, \dots, z_{n}) \thickapprox 0$.

It only remains to prove (iv)$\Rightarrow$(i). Suppose that the equation $p(z_{1}, \dots, z_{n}) \thickapprox 0$ admits 
an integer solution. Recall that the free commutative ring with free generators $\{z_{1}, \dots, z_{n}\}$ is the polynomial ring $\mathds{Z}[z_{1}, \dots, z_{n}]$. Since $\mathds{Z}$ is a subalgebra of 
$\mathds{Z}[z_{1}, \dots, z_{n}]$, this implies that there are constant terms (in the language of commutative rings) 
$\alpha_{1}, \dots, \alpha_{n}$ such that $\mathds{Z}[z_{1}, \dots, z_{n}]\vDash p(\alpha_{1}, \dots, \alpha_{n}) \thickapprox 0$ and, therefore, that $\class{CR} \vDash p(\alpha_{1}, \dots, \alpha_{n}) \thickapprox 0$. From point 3 of Lemma \ref{DeductionLemma} it follows that $p(\alpha_{1}, \dots, \alpha_{n}) \sineq_{\class{CR}} 0$. By Theorem \ref{Finite} this is equivalent to the fact that $p(\alpha_{1}, \dots, \alpha_{n}) \sineq_{\mathcal{CR}} 0$. Therefore we can apply Lemma \ref{FiniteProofs} obtaining a finite sequence $\langle \gamma_{1}, \gamma_{2}, \dots, \gamma_{m}\rangle$, where $\gamma_{1}= p(\alpha_{1}, \dots, \alpha_{n})$, $\gamma_{m} = 0$ and for every $k < m$ there is a rule $\alpha \sineq\beta \in \mathcal{CR}$ and a substitution $\sigma$ (in the language of commutative rings) such that $\{\sigma\alpha, \sigma \beta\} = \{\gamma_{k}, \gamma_{k+1}\}$. Observe that $\sigma$ can be regarded as a substitution in the language of $\LL(p)$ as well. Moreover, recall that $\alpha\leftrightarrow \beta$ is an axiom of $\LL(\? p )$ by (CR). Therefore, by structurality, we obtain $\emptyset \vdash_{\LL(\? p )}\sigma\alpha \leftrightarrow \sigma\beta$. Applying (S) if necessary, this yields that $\emptyset \vdash_{\LL(\? p )} \gamma_{k} \leftrightarrow \gamma_{k+1}$. Hence we proved that $\emptyset \vdash_{\LL(\? p )} \gamma_{k} \leftrightarrow \gamma_{k+1}$ for every $k \leq m-1$. Applying $m-1$ times (T) we conclude that
\begin{equation}\label{OurNiceTheorem}
\emptyset \vdash_{\LL(\? p )} p(\alpha_{1}, \dots, \alpha_{n}) \leftrightarrow 0.
\end{equation}
Recall that the variables $x$ and $y$ do not appear in the equation $p(z_{1}, \dots, z_{n}) \thickapprox 0$ therefore, we can safely consider the substitution $\sigma$ that sends $z_{i}$ to $\alpha_{i}$ for every $i \leq n$ and leaves the other variables untouched. By (MP'), (A3'), (G') and  (\ref{OurNiceTheorem}) we obtain that  
\[
x, x\leftrightarrow y \vdash_{\LL(\? p )}y \quad x \sineq_{\LL(\? p )} x \leftrightarrow ( x \leftrightarrow x) \text{ and } x, y \vdash_{\LL(\? p )} x \leftrightarrow y.
\]
Now it is easy to see that what we have proved is just the syntactic characterization of algebraizability  of Theorem \ref{SynAlg} for $\btau(x) \coloneqq \{ x \thickapprox x \leftrightarrow x\}$ and $\brho(x, y) \coloneqq \{ x \leftrightarrow y \}$. Therefore, with an application of Lemma \ref{AssertionalDelta}, we conclude that $\LL(\? p )$ is finitely regularly algebraizable.
\end{proof}
\begin{Corollary}\label{Cor : Consistency}
The logic $\LL(p)$ is consistent for every Diophantine equation $p \thickapprox 0$.
\end{Corollary}
\begin{proof}
It is easy to see that the matrix $\langle \boldsymbol{Z}, \{ s \} \rangle$ defined in the above proof is a model of $\LL(p)$. Since $\{ s \} \ne Z$, we conclude that $\LL(p)$ is consistent.
\end{proof}

Observe that every class of the Leibniz hierarchy is contained either into the class of protoalgebraic logic or into the one of logics $\LL$ whose truth sets are implicitly definable in $\ModstarL$. Moreover every class of the Leibniz hierarchy contains the one of finitely regularly algebraizable logics. Keeping this is mind, Theorem \ref{Thm:LeibnizUndecidable} shows that the problem of determining whether a logic of the form $\LL(p)$ belong to a \textit{given} level of the Leibniz hierarchy is equivalent to the one of determining whether it belongs to \text{any} level of the Leibniz hierarchy. This does not contradicts the fact that \textit{in general} these two problems are different. By means of this peculiar feature of logics of the form $\LL(p)$, we are able to establish at once the undecidability of the various problems (one for each level of the Leibniz hierarchy) of determining whether a logic belong to a given level of the Leibniz hierarchy.

\begin{Theorem}\label{UndecAlg} Let $\class{K}$ be a level of the Leibniz hierarchy in Figure \ref{Hierarchy}. The problem of determining whether the logic of a given consistent finite Hilbert calculus in a finite language belongs to $\class{K}$ is undecidable.
\end{Theorem}
\begin{proof}
Suppose towards a contradiction that there is an algorithm $\class{A}_{1}$ that, given a consistent finite Hilbert calculus in a finite language, determines whether its logic belongs to $\class{K}$. Then we define a new algorithm $\class{A}_{2}$ as follows: given a Diophantine equation $p \thickapprox 0$, we construct the logic $\LL(p)$ and check with $\class{A}_{1}$ whether it belongs to $\class{K}$. Observe that we can do this, since $\LL(p)$ is consistent by Corollary \ref{Cor : Consistency}. Then in the positive case $\class{A}_{2}$ returns \textit{yes}, while \textit{no} otherwise. Observe that $\class{K}$ contains the class of finitely regularly algebraizable logics and is included either in the class of protoalgebraic logics or in the class of logics $\LL$ whose truth sets are implicitly definable in $\ModstarL$. Therefore we can apply Theorem \ref{Thm:LeibnizUndecidable}, yielding that:
\[
\LL(p) \in \class{K} \Longleftrightarrow p \thickapprox 0 \text{ has an integer solution}.
\]
Therefore $\class{A}_{2}$ would provide a decision procedure for Hilbert's tenth problem. Since we know that such a procedure does not exist, we obtain a contradiction as desired.
\end{proof}

\section{The classification problem in the Frege hierarchy}

Now we move our attention to the problem of classifying logics in the Frege hierarchy, which deals with several kinds of replacement properties. 
Our goal is to prove a result analogous to the one obtained for the Leibniz hierarchy. 
More precisely, we will show that the problem of classifying the logic of a consistent finite Hilbert 
calculus in a finite language in the Frege hierarchy is in general undecidable (Theorem \ref{Thm:FregeUndecidable}). Remarkably, our argument shows that this classification problem remains undecidable even if we restrict our attention to Hilbert calculi that determine a finitary algebraizable logic. Our strategy is as follows. Consider a finitely based variety $\class{V}$ of finite type, whose equational theory in one variable is undecidable.\footnote{Observe that examples of finitely based varieties of finite type, whose equational theory in one variable is undecidable, are known. For instance, the variety of \textit{relation algebras} is of this kind, as shown in \cite[Section~8.5(viii)]{TaGi87}. Other examples are described in \cite{AIMa66}.} We will reduce the problem of determining the equational theory of $\class{V}$ in a single variable to the problem of classifying logics of finite consistent Hilbert calculi in the Frege hierarchy.

 Let $\Phi$ be a finite equational basis for the variety $\class{V}$, and let $Eq(x)$ be the set of all equations in the language of $\class{V}$ in  a single variable $x$. We will introduce a new logic $\LL(\alpha, \beta)$ for each equation $\alpha \thickapprox \beta \in Eq(x)$. 

\begin{Definition}
Let $\alpha \thickapprox \beta \in Eq(x)$ and let $\mathscr{L}$ be the language of $\class{V}$, expanded with two new connectives $\Box$ and $\to$, respectively unary and binary. Then $\LL(\alpha, \beta)$ is the logic in the language $\mathscr{L}$ axiomatized by the following Hilbert calculus:
\begin{align*}
\emptyset \vdash & x \to x \tag{\textup{R}}\\
x, x \to y \vdash & y \tag{\textup{MP}}\\
\{ x_{i} \to y_{i}, y_{i} \to x_{i} : 1 \leq i \leq n\}\vdash & f(x_{1}, \dots, x_{n}) \to f(y_{1}, \dots, y_{n})\tag{Rep}\\
x \sineq  & \Box x \to x, x \to \Box x \tag{A3}\\
\emptyset \vdash & \epsilon \to \delta, \delta \to \epsilon \tag{\textup{V}}\\
\alpha( \varphi ) \to \beta( \varphi ) \vdash & \varphi \tag{W}
\end{align*}
for every $\epsilon \thickapprox \delta \in \Phi$, every $n$-ary connective $f$, and every formula $\varphi$ of the following list:
\begin{align*}
\varphi_{1} \coloneqq & x \to (y \to x)\\
\varphi_{2}  \coloneqq & (x \to (y \to z)) \to ( ( x \to y) \to ( x \to z))\\
\varphi_{3}  \coloneqq & (x_{1} \to y_{1}) \to (  (y_{1} \to x_{1}) \to    (\dots \to((x_{n} \to y_{n}) \to\\
& \to ((y_{n} \to x_{n}) \to (f(\vec{x}) \to f(\vec{y}))))\dots))\\
\varphi_{4}  \coloneqq & x \to ( x \to \Box x)\\
\varphi_{5}  \coloneqq & x \to ( \Box x \to x)\\
\varphi_{6}  \coloneqq & (\Box x \to x) \to ( (x \to \Box x ) \to x)
\end{align*}
again with $f$ ranging over all $n$-ary for every $n \in \omega$.
\end{Definition}

Observe that the syntactic characterization of algebraizability (Theorem \ref{SynAlg}) implies that $\LL(\alpha, \beta)$ is finitely algebraizable through
\[
\btau(x) = \{ x \thickapprox \Box x \} \text{ and }\brho(x, y) = \{ x \to y, y \to x \}.
\]
The next result expresses the relation between the validity of the equation $\alpha \thickapprox \beta$ in $\class{V}$ and the location of the logic $\LL(\alpha, \beta)$ in the Frege hierarchy.

\begin{Theorem}\label{Thm:FregeUndecidable}
Let $\alpha \thickapprox \beta \in Eq(x)$. The following conditions are equivalent:
\benroman
\item $\LL(\alpha, \beta)$ is fully Fregean.
\item $\LL(\alpha, \beta)$ is selfextensional.
\item $\class{V} \vDash \alpha \thickapprox \beta$.
\eroman
\end{Theorem}

\begin{proof}
(i)$\Rightarrow$(ii) is straightforward. Then we turn to prove (ii)$\Rightarrow$(iii). Suppose towards a contradiction that $\LL(\alpha, \beta)$ is selfextensional and that $\class{V} \nvDash \alpha \thickapprox \beta$. Then let $\A \in \class{V}$ be the free algebra with countably many free generators. It is well known that the universe of $\A$ consists of congruence classes of terms, that we denote by $\llbracket \varphi \rrbracket, \llbracket \psi \rrbracket\dots$ In particular, we denote the free generators of $\A$  by $\{ \llbracket z \rrbracket : z \in \Var \}$. In particular, consider the three distinct elements  $\llbracket x \rrbracket, \llbracket y \rrbracket$ and $\llbracket v \rrbracket$ of $\A$ corresponding to the three distinct variables $x, y$ and $v$. Then  expand $\A$ with two new unary and binary operations $\Box$ and $\to$ defined as
\[
\Box a \coloneqq \left\{ \begin{array}{ll}
 \llbracket y \rrbracket & \text{if $a = \llbracket x \rrbracket$}\\
 a & \text{otherwise}
  \end{array} \right.    
\quad 
a \to b \coloneqq \left\{ \begin{array}{ll}
 \llbracket x \rrbracket & \text{if $a \ne b$}\\
 \llbracket y \rrbracket & \text{if $a = b  \ne \llbracket x \rrbracket$}\\
 \llbracket v \rrbracket & \text{if $a = b = \llbracket x \rrbracket$}
  \end{array} \right.  
\]
for every $a, b \in A$. Let $\B$ be the result of the expansion and put $F \coloneqq B \smallsetminus \{ \llbracket x \rrbracket \}$.

\begin{claim}\label{Claim:1a}
$\langle \B, F \rangle \in \ModstarL(\alpha, \beta)$.
\end{claim}

First observe that $\langle \B, F \rangle$ is a reduced matrix. To prove this, consider two different elements $a, b \in B$ and the polynomial function $p(x) \coloneqq x \to a$. It is easy to see that $p(a) \in F$ and $p(b) \notin F$. By point 1 of Lemma \ref{Polynomial} we conclude that $\langle a, b \rangle \notin \leibniz^{\B}F$ and, therefore, that $\langle \B, F \rangle$ is reduced. Then we turn to prove that it is a model of $\LL(\alpha, \beta)$. The closure of $F$ under the rules axiomatizing $\LL(\alpha, \beta)$ is easily proved, except perhaps for the set of rules (W). Consider a rule $\alpha(\varphi) \to \beta( \varphi) \vdash \varphi$ in (W) and a homomorphism $h\colon \Fm \to \B$ such that $h \alpha(\varphi) \to h \beta(\varphi) \in F$. Clearly we have that $h\alpha(\varphi) = h\beta(\varphi)$. Since $\class{V} \nvDash \alpha \thickapprox \beta$ and $\A$ is the free relation algebra, we know that $h(\varphi) \notin  \{ \llbracket z \rrbracket : z \in \Var \}$. Together with the definition of $F$, this implies that $h(\varphi) \in F$.

Now observe that both $x \to x$ and $y \to y$ are 
instances of (\textup{R}) and, therefore, they are theorems of $\LL(\alpha, \beta)$. In particular this means 
that $x \to x \sineq_{\LL(\alpha, \beta)} y \to y$. Since $\LL(\alpha, \beta)$ 
is selfextensional, we obtain that $\AlgL(\alpha, \beta) \vDash x \to x \thickapprox y \to y$ by Lemma \ref{EquationalFrege}. Let $g\colon \Fm \to \B$ the the natural surjection, which sends $z$ to $\llbracket z \rrbracket$ for every $z \in Var$. It is easy to see that
\[
g(x \to x) = \llbracket v \rrbracket \ne  \llbracket y \rrbracket = g(y \to y).
\]
But from Claim \ref{Claim:1a} it follows that $\B \in \AlgL$. Therefore we reach a contradiction as desired.

It only remains to prove part (iii)$\Rightarrow$(i). Suppose that $\class{V} \vDash \alpha \thickapprox \beta$. Recall that $\LL( \alpha, \beta)$ is algebraizable with set of congruence formulas $\brho(x, y) = \{ x \to y, y \to x \}$. By Lemma \ref{Lem:MovingEqIntoFm} we have that
\[
\emptyset \vdash_{\LL(\alpha, \beta)} \brho(\gamma, \eta) \Longleftrightarrow \AlgL(\alpha, \beta) \vDash \gamma \thickapprox \eta
\]
for every $\gamma, \eta \in Fm$. Applying this observation to the axiom (\textup{V}), we conclude 
that $\AlgL(\alpha, \beta)$ is a class of expansions of algebras in $\class{V}$. Together with the fact that $\class{V} \vDash \alpha \thickapprox \beta$, this implies that $\alpha \to \beta$ is a theorem of $\LL(\alpha, \beta)$. Then consider any $i \leq 7$. By structurality $\alpha(\varphi_{i}) \to \beta(\varphi_{i})$ is a theorem of $\LL(\alpha, \beta)$. Therefore we can apply the rules (W), obtaining that also $\varphi_{i}$ is a theorem of $\LL(\alpha, \beta)$. 

Now, let $\LL$ be the logic axiomatized only by $\{ \varphi_{i} : i \leq 6 \}$, (V) and (MP).

\begin{claim}\label{Claim:1b}
 $\LL = \LL(\alpha, \beta)$.
\end{claim}

Clearly we have that $\LL \leq \LL( \alpha, \beta)$.  To show the other direction, observe that $\LL$ satisfies the rules (W). Moreover, with an application of (MP) and $\varphi_{3}, \varphi_{4}, \varphi_{5}, \varphi_{6}$, it is easy to see that $\LL$ satisfies (Rep) and (A3) too. It only remains to show that $\LL$ satisfies (R). But observe that $\LL$ is an expansion of the implication fragment $\mathcal{IPC}_{\to}$ of propositional intuitionistic logic, which is axiomatized by $\varphi_{1}, \varphi_{2}$ and (MP) \cite[Example 5.91]{AAL-AIT-f}. Since $x \to x$ is a theorem of $\mathcal{IPC}_{\to}$, we conclude that $\LL$ satisfies (R) and, therefore, that $\LL( \alpha, \beta) \leq \LL$. Thus we established the claim.

It is well known that a finitary logic in a language containing $\to$ satisfies the classical version of the deduction-detachment theorem if and only if it is an axiomatic extension of the logic defined in its language by the axioms $\varphi_{1}$ and $\varphi_{2}$ and the rule (MP) \cite[Theorem 2.4.2]{W88}. Together with Claim \ref{Claim:1b}, this implies that
\[
\Gamma \vdash_{\LL( \alpha, \beta)} \gamma \to \psi \Longleftrightarrow \Gamma, \gamma \vdash_{\LL( \alpha, \beta)} \psi
\]
for every $\Gamma \cup \{ \gamma, \psi \} \subseteq Fm$. In particular, this means that for every $\gamma, \psi \in Fm$:
\[
\gamma\sineq_{\LL( \alpha, \beta)} \psi \Longleftrightarrow \emptyset \vdash_{\LL( \alpha, \beta)} \brho( \gamma, \psi) \Longleftrightarrow \AlgL \vDash \gamma \thickapprox \psi.
\]
This fact easily implies that $\sineq_{\LL( \alpha, \beta)}$ is a congruence of $\Fm$. Hence we conclude that $\LL( \alpha, \beta)$ is selfextensional. The fact that $\LL( \alpha, \beta)$ is also fully Fregean, is a consequence of two general results of abstract algebraic logic. First, every finitary selfextensional logic with a classical deduction-detachment theorem is fully selfextensional \cite[Theorem 4.46]{FJa09} (see also \cite{Ja05}).\footnote{In these references the classical deduction-detachment theorem is called \textit{uniterm} deduction-detachment theorem. The expression \textit{uniterm} refers to the fact that this deduction-detachment theorem is witnessed by a single formula in two variables in contrast to the cases where the deduction-detachment theorem is witnesses by a set of formulas as in \cite{BP-AAL-DDT,CzePi04}.} Second, every truth-equational fully selfextensional logic is also fully Fregean \cite[Theorem 22]{AFRM15}.
\end{proof}

We are now ready to state the main result on the classification of syntactically presented logics in the Frege hierarchy.

\begin{Theorem}\label{Thm:FregeUndecidable2} Let $\class{K}$ be a level of the Frege hierarchy. The problem of determining whether the logic of a given finite consistent Hilbert calculus in a finite language belongs to $\class{K}$ is undecidable. Moreover, the problem remains undecidable when restricted to the classification of finite consistent Hilbert calculi that determine a finitely algebraizable logic.
\end{Theorem}

\begin{proof}
Observe that the logic $\LL(\alpha, \beta)$ was defined starting from an arbitrary variety $\class{V}$ of finite type, which is finitely based and whose equational theory in one variable is undecidable. We can assume w.l.o.g.\ that $\class{V}$ is the variety $\class{RA}$ of relation algebras. Under this assumption, we claim that the logic $\LL(\alpha, \beta)$ is consistent for every $\alpha \thickapprox \beta \in Eq(x)$.  The fact that $\LL(\alpha, \beta)$ is consistent when $\class{RA} \nvDash \alpha \thickapprox \beta$ has been proved in part (ii)$\Rightarrow$(iii) of Theorem \ref{Thm:FregeUndecidable}. Then consider the case where $\class{RA} \vDash \alpha \thickapprox \beta$. In part (iii)$\Rightarrow$(i) of Theorem \ref{Thm:FregeUndecidable}, we showed that $\LL( \alpha, \beta)$ is the logic axiomatized by only by $\{ \varphi_{i} : i \leq 6 \}$, (V) and (MP). Then consider any non-trivial Boolean algebra $\A = \langle A, \land, \lor, \lnot\rangle$ and expand it to a relation algebra $\langle A, \land, \lor, \lnot, \cdot, \?^{\smallsmile}\!, 1\rangle$ by interpreting $\cdot$ as $\land$, $^{\smallsmile}$ as the identity map, and $1$ as the top element. Moreover, let $\to$ be the usual Boolean implication and $\Box$ be the constant function with value $1$. Let $\B$ be the result of the expansion. It is very easy to see that $\langle \B, \{ 1 \} \rangle$ is a model of $\{ \varphi_{i} : i \leq 6 \}$, (V) and (MP). Hence $\LL(\alpha, \beta)$ is consistent. This establishes our claim.

Now, we turn back to the main proof, again assuming that $\class{V} = \class{RA}$. Let $\class{K}$ be a level of the Frege hierarchy. Suppose towards a contradiction that there is an algorithm $\class{A}_{1}$ which, given a finite consistent Hilbert calculus in a finite language that moreover determines a finitely algebraizable logic, decides whether its logic belongs to $\class{K}$. Then we define a new algorithm $\class{A}_{2}$ as follows: given an equation $\alpha \thickapprox \beta \in Eq(x)$, we construct the logic $\LL(\alpha, \beta)$ and check with $\class{A}_{1}$ if it belongs to $\class{K}$. In the positive case $\class{A}_{2}$ returns \textit{yes}, while \textit{no} otherwise. Observe that we can do this, since $\LL(\alpha, \beta)$ is finitely algebraizable and, by the claim,  consistent.  Since $\class{K}$ contains the class of selfextensional logics and is included in the of class of fully Fregean ones, we can apply Theorem \ref{Thm:FregeUndecidable} obtaining that
\[
\LL(\alpha, \beta) \in \class{K} \Longleftrightarrow \class{RA} \vDash \alpha \thickapprox \beta.
\]
Therefore $\class{A}_{2}$ would provide a decision procedure for the validity in $\class{RA}$ of equations in one variable. But this contradicts the known fact that such a procedure does not exist \cite[Section~8.5(viii)]{TaGi87}.
\end{proof}

\section*{Appendix}

\begin{proof}[Proof of Lemma 4.5]
First we check that $\sineq_{\mathcal{CR}}$ preserves the connective $-$. To do this, consider $\psi, \varphi \in Fm$ such that $\psi \sineq_{\mathcal{CR}} \varphi$: we have to prove that $- \psi \sineq_{\mathcal{CR}} - \varphi$. By Lemma \ref{FiniteProofs} it will be enough to check that $- \alpha \sineq_{\mathcal{CR}} - \beta$ for every $\alpha \sineq \beta \in \mathcal{CR}$. We shall make use of some deductions. In particular, by suitable substitutions and applying (N) and (O) to (B), (D), (E), (H), (I), (L) and (M), we obtain that
\begin{align*}
x \cdot y&\sineq_{\mathcal{CR}} y \cdot x \tag{\textup{B'}}\\
(x + y) + z &\sineq_{\mathcal{CR}} x + ( y + z ) \tag{\textup{D'}}\\
x + y  &\sineq_{\mathcal{CR}}  y + x \tag{\textup{E'}}\\
 x \cdot ( y + z ) &\sineq_{\mathcal{CR}}  ( x \cdot y ) + (x \cdot z ) \tag{\textup{H'}}\\
-(x + y ) &\sineq_{\mathcal{CR}} - x + - y \tag{\textup{I'}}\\
w + -(x \cdot y ) &\sineq_{\mathcal{CR}}w + ( - x \cdot  y) \tag{\textup{L'}}\\
w + -(x \cdot y ) &\sineq_{\mathcal{CR}}w + (  x \cdot - y) \tag{\textup{M'}}
\end{align*}
Then suppose that we are given a rule of our calculus (X), different from (N) and (O). Then (X) is of the form $w + (u \cdot \epsilon) \sineq w + ( u \cdot \delta)$ for some formulas $\epsilon$ and $\delta$. But we have that
\begin{align*}
- (w + (u \cdot \epsilon)) &\sineq_{\mathcal{CR}} -w + -(u \cdot \epsilon) \tag{\textup{I'}}\\
&\sineq_{\mathcal{CR}} -w + (-u \cdot \epsilon) \tag{\textup{L'}}\\
&\sineq_{\mathcal{CR}} -w + (-u \cdot \delta) \tag{\textup{X}}\\
&\sineq_{\mathcal{CR}} -w + -(u \cdot \delta) \tag{\textup{L'}}\\
&\sineq_{\mathcal{CR}} - (w + (u \cdot \delta)) \tag{\textup{I'}}
\end{align*}
Therefore it only remains to prove the cases of (N) and (O). This is what we do now:
\begin{align*}
- (0 + x) &\sineq_{\mathcal{CR}} 0 + - (0 + x) \tag{\textup{N}}\\
&\sineq_{\mathcal{CR}} 0 + (1 \cdot - (0 +x)) \tag{\textup{O}}\\
&\sineq_{\mathcal{CR}} 0 + - (1 \cdot (0+x)) \tag{\textup{M'}}\\
&\sineq_{\mathcal{CR}} 0 + (-1 \cdot (0 + x)) \tag{\textup{L'}}\\
&\sineq_{\mathcal{CR}}0 + (-1\cdot (x + 0))\tag{\textup{E}}\\
&\sineq_{\mathcal{CR}} 0 + (-1 \cdot x) \tag{\textup{F}}\\
&\sineq_{\mathcal{CR}} 0 + -(1 \cdot x)  \tag{\textup{L'}}\\
&\sineq_{\mathcal{CR}} 0 + (1 \cdot - x) \tag{\textup{M'}}\\
&\sineq_{\mathcal{CR}} 0 + - x \tag{\textup{O}}\\
&\sineq_{\mathcal{CR}} - x \tag{\textup{N}}
\end{align*}
and
\begin{align*}
- (x + ( 1 \cdot y)) &\sineq_{\mathcal{CR}} - x + - (1 \cdot y) \tag{\textup{I'}}\\
&\sineq_{\mathcal{CR}} - x + (1 \cdot - y) \tag{\textup{M'}}\\
&\sineq_{\mathcal{CR}} - x + - y \tag{\textup{0}}\\
&\sineq_{\mathcal{CR}} - ( x + y) \tag{\textup{I'}}\\
\end{align*}
Therefore we conclude that $\sineq_{\mathcal{CR}}$ preserves $-$.

Then we turn to prove that $\sineq_{\mathcal{CR}}$ preserves $+$ and $\cdot$ too. In order to do this, it will be enough to show that if $\psi \sineq_{\mathcal{CR}} \varphi$, then $\chi + \psi \sineq_{\mathcal{CR}} \chi + \varphi$ and $\chi \cdot \psi \sineq_{\mathcal{CR}} \chi \cdot \varphi$ for every formula $\chi$. Let us explain briefly why. Suppose that this condition, call it (Y), holds. Then consider $\psi_{1}, \psi_{2}, \varphi_{1}, \varphi_{2} \in Fm$ such that $\psi_{1} \sineq_{\mathcal{CR}} \varphi_{1}$ and $\psi_{2} \sineq_{\mathcal{CR}} \varphi_{2}$. We would have that
\begin{align*}
\psi_{1} + \psi_{2}&\sineq_{\mathcal{CR}}\psi_{1} + \varphi_{2} \tag{\textup{Y}}\\
&\sineq_{\mathcal{CR}} \varphi_{2} + \psi_{1} \tag{\textup{E'}}\\
&\sineq_{\mathcal{CR}} \varphi_{2} + \varphi_{1} \tag{\textup{Y}}\\
&\sineq_{\mathcal{CR}} \varphi_{1} + \varphi_{2} \tag{\textup{E'}}\\
\end{align*}
and
\begin{align*}
\psi_{1} \cdot \psi_{2}&\sineq_{\mathcal{CR}} \psi_{1} \cdot \varphi_{2} \tag{\textup{Y}}\\
&\sineq_{\mathcal{CR}} \varphi_{2} \cdot \psi_{1} \tag{\textup{B'}}\\
&\sineq_{\mathcal{CR}} \varphi_{2} \cdot \varphi_{1} \tag{\textup{Y}}\\
&\sineq_{\mathcal{CR}} \varphi_{1} \cdot \varphi_{2} \tag{\textup{B'}}\\
\end{align*}
concluding the proof. Therefore we turn to prove that if $\psi \sineq_{\mathcal{CR}} \varphi$, then $\chi + \psi \sineq_{\mathcal{CR}} \chi + \varphi$ and $\chi \cdot \psi \sineq_{\mathcal{CR}} \chi \cdot \varphi$ for every formula $\chi$. Suppose that $\psi \sineq_{\mathcal{CR}} \varphi$ and consider an arbitrary formula $\chi$. By Lemma \ref{FiniteProofs}, to prove that $\chi + \psi \sineq_{\mathcal{CR}} \chi + \varphi$ and $\chi \cdot \psi \sineq_{\mathcal{CR}} \chi \cdot \varphi$, it will be enough to check that $\chi + \alpha \sineq_{\mathcal{CR}} \chi + \beta$ and $\chi \cdot \alpha \sineq_{\mathcal{CR}} \chi \cdot \beta$ for every rule $\alpha \sineq \beta$ in $\mathcal{CR}$. Then suppose that we are given a rule (X) of $\mathcal{CR}$, different from (N) and (O). Then (X) is of the form $w + (u \cdot \epsilon) \sineq w + ( u \cdot \delta)$ for some formulas $\epsilon$ and $\delta$. We have that
\begin{align*}
\chi + (w + (u \cdot \epsilon)) &\sineq_{\mathcal{CR}} (\chi + w) + (u \cdot \epsilon) \tag{\textup{D'}}\\
&\sineq_{\mathcal{CR}} (\chi + w) + (u \cdot \delta) \tag{\textup{X}}\\
&\sineq_{\mathcal{CR}} \chi + (w + (u \cdot \delta)) \tag{\textup{D'}}
\end{align*}
and
\begin{align*}
\chi \cdot (w + (u \cdot \epsilon)) &\sineq_{\mathcal{CR}} (\chi \cdot w) + (\chi \cdot (u \cdot \epsilon)) \tag{\textup{H'}}\\
&\sineq_{\mathcal{CR}} (\chi \cdot w) + (1 \cdot(\chi \cdot (u \cdot \epsilon))) \tag{\textup{O}}\\
&\sineq_{\mathcal{CR}} (\chi \cdot w) + (1 \cdot((\chi \cdot u) \cdot \epsilon)) \tag{\textup{A}}\\
&\sineq_{\mathcal{CR}} (\chi \cdot w) + ((\chi \cdot u) \cdot \epsilon) \tag{\textup{O}}\\
&\sineq_{\mathcal{CR}} (\chi \cdot w) + ((\chi \cdot u) \cdot \delta) \tag{\textup{X}}\\
&\sineq_{\mathcal{CR}} (\chi \cdot w) + (1 \cdot((\chi \cdot u) \cdot \delta)) \tag{\textup{O}}\\
&\sineq_{\mathcal{CR}} (\chi \cdot w) + (1 \cdot(\chi \cdot (u \cdot \delta))) \tag{\textup{A}}\\
&\sineq_{\mathcal{CR}} (\chi \cdot w) + (\chi \cdot (u \cdot \delta)) \tag{\textup{O}}\\
&\sineq_{\mathcal{CR}} \chi \cdot (w +  (u \cdot \delta)) \tag{\textup{H'}}
\end{align*}
Therefore it only remains to check cases (N) and (O). For what concerns (N) we have that:
\begin{align*}
\chi + (0 + x) &\sineq_{\mathcal{CR}} (0 + x) + \chi \tag{\textup{E'}}\\
&\sineq_{\mathcal{CR}} 0 + ( x + \chi)  \tag{\textup{D'}}\\
&\sineq_{\mathcal{CR}} x + \chi \tag{\textup{N}}\\
&\sineq_{\mathcal{CR}} \chi + x \tag{\textup{E'}}
\end{align*}
and
\begin{align*}
\chi \cdot (0 + x) &\sineq_{\mathcal{CR}} 0 + (\chi \cdot (0 + x)) \tag{\textup{N}}\\
&\sineq_{\mathcal{CR}} 0 + (\chi \cdot ( x + 0 )) \tag{\textup{E}}\\
&\sineq_{\mathcal{CR}} 0 + (\chi \cdot  x ) \tag{\textup{F}}\\
&\sineq_{\mathcal{CR}} \chi \cdot  x  \tag{\textup{N}}
\end{align*}
Then we turn to prove the case of (O). We have that
\begin{align*}
\chi + (x + (1\cdot y)) &\sineq_{\mathcal{CR}} (\chi + x) + (1 \cdot y) \tag{\textup{D'}}\\
&\sineq_{\mathcal{CR}} (\chi + x) + y  \tag{\textup{O}}\\
&\sineq_{\mathcal{CR}} \chi + ( x + y) \tag{\textup{D'}}
\end{align*}
and
\begin{align*}
\chi \cdot (x + (1\cdot y)) &\sineq_{\mathcal{CR}} (\chi \cdot x) + ( \chi \cdot (1 \cdot y)) \tag{\textup{H'}}\\
&\sineq_{\mathcal{CR}} (\chi \cdot x) + ( \chi \cdot (y \cdot 1))  \tag{\textup{B}}\\
&\sineq_{\mathcal{CR}} (\chi \cdot x) + ( \chi \cdot y) \tag{\textup{C}}\\
&\sineq_{\mathcal{CR}} \chi \cdot (x +  y) \tag{\textup{H'}}
\end{align*}
This concludes the proof that $\mathcal{CR}$ is selfextensional.
\end{proof}

\

\paragraph{\bfseries Acknowledgements}
The idea of considering the problem of classifying logics in the Leibniz and Frege hierarchies from a computational point of view originates from a list of open problems in abstract algebraic logic presented by F\'elix Bou in the Seminar on Non-Classical Logics of the University of Barcelona. I am very grateful to him for bringing my attention of this topic. Thanks are due also to Professors Josep Maria Font, Jos\'e Gil-F\'erez, Ramon Jansana and James Raftery for their useful comments on early versions of the paper, which contributed to improve a lot its presentation, and to Eduardo Hermo for discussing its topic when the work was at a preliminary stage. Finally, I wish to thank the anonymous referee, whose comments and suggestions improved further the presentation of the paper. The author was supported by the project GA13-14654S of the Czech Science Foundation.

\end{document}